\documentclass[12pt,letterpaper]{amsart}
\usepackage[english]{babel}
\usepackage[dvips]{graphicx}
\usepackage{amscd}
\usepackage{latexsym}
\usepackage{amsfonts}
\newtheorem{theorem}{Theorem}
\newtheorem{proposition}[theorem]{Proposition}
\newtheorem{lemma}[theorem]{Lemma}
\newtheorem{corollary}[theorem]{Corollary}
\theoremstyle{definition}

\begin{document}
\title{The Strong Symmetric Genus of the Hyperoctahedral Groups}
\author{Michael A. Jackson}
\maketitle
\section{Introduction}
In the study of Reimannian manifolds, it is natural to consider the finite groups, which act as automorphisms 
of the manifold. On the other hand, given a finite group, one may consider the topological surfaces on which the 
group acts faithfully as a group of automorphisms. There are several natural invariants assigned 
to a group that are associated to the action of that group on compact orientable surfaces. The first of these 
invariants is the genus $\gamma (G)$ of the group $G$, which is the smallest genus of a surface on which 
some Cayley graph for $G$ can be embedded. The symmetric genus $\sigma (G)$ of a finite group 
$G$ is the smallest genus of a surface on which $G$ acts faithfully as a group of automorphisms. 
The strong symmetric genus $\sigma ^{0} (G)$ is the smallest genus of a surface on which $G$ acts 
faithfully as a group of orientation preserving automorphisms. (See \cite{gt:tgt} chapter 6.) It is clear 
that these invariants satisfy the inequality $\gamma (G) \leq \sigma (G) \leq \sigma ^{0} (G)$.

Perhaps the most classical of these parameters is the strong symmetric genus. It was first 
considered by Burnside \cite{b:tgfo}, who included the further restriction that the quotient space under the 
action be a sphere. The strong symmetric genus is what we will concern ourselves 
with in this paper. Various results have been shown involving the strong symmetric genus. 
For example, all groups $G$ such that $\sigma ^{0} (G)\leq 3$ have been computed (see 
\cite{b:cfga,mz:gsssg}). Also the strong symmetric genus of several infinite families of finite groups 
have been found, such as the alternating and symmetric groups \cite{c:gasg,c:mgasg,c:srqtg}, the 
groups $PSL_2 (q)$ \cite{gs:gop,gs:rpsl}, and the groups $SL_2(q)$ \cite{v:gos}. In addition, the 
strong symmetric genus has been found for the sporadic finite simple groups \cite{cww:sgsg,w:mhg,w:sgfg,w:sgbm}.

Since the strong symmetric genus of the alternating and symmetric groups is known, it is natural to 
try to find the strong symmetric genus for other infinite families of finite groups. An example of such an infinite 
family is the collection of hyperoctahedral groups. Recall that the hyperoctahedral groups are the finite Coxeter 
groups of type $B_n$. The hyperoctahedral group $B_n$ for $n\geq 3$ is defined as the group of symmetries of the 
$n$-dimensional cube. In this paper we will prove the following theorem and a corollary about the hyperoctahedral groups:

\begin{theorem}\label{thm:main}
For all $n\geq 3$, except $n=$5, 6, and 8, the hyperoctahedral group $B_n$ is a quotient of the 
triangle group $T(2,4,6)=\langle x,y,z | x^2=y^4=z^6=xyz=1\rangle$. 
In the exceptional cases, $B_5$, $B_6$ and $B_8$ are quotients of the triangle groups $T(2,4,10)$, 
$T(2,6,6)$, and $T(2,4,8)$ respectively.
\end{theorem}

\begin{corollary}\label{cor:main}
For all $n\geq 3$, except $n=$5, 6, and 8, the strong symmetric genus of the hyperoctahedral group $B_n$ is 
$\frac{n! 2^n}{24}+1=\frac{n! 2^{n-3}}{3}+1$. The strong symmetric genus in the exceptional cases $B_5$, $B_6$, 
and $B_8$ are 289, 3841, and 645,121 respectively.
\end{corollary}

The investigation of the strong symmetric genus of the hyperoctahedral groups 
grew out of a VIGRE Working Group at The Ohio State University during 
the fall of 2002. In this working group, we were looking to find the strong symmetric genus of some finite 
Coxeter groups. The exceptional cases of Theorem \ref{thm:main} were found by the working group \cite{vwg:ggt}.

\section{Subgroups of the Hyperoctahedral groups}
For the present time we will fix an $n \geq 3$ and look at the hyperoctahedral group $B_n$. 
As we have previously said, $B_n$ is defined as the group of symmetries of the $n$-dimensional 
cube. This group is also described as the wreath product $\mathbb{Z}_2 \wr \Sigma _n$. 
Another description of this group is as the group of all $n\times n$ signed permutation matrices. For our 
purposes, we will use the wreath product description as adopted by V. S. Sikora \cite{s:tbhg}. For an 
element of $B_n$, we will write a tuple $[\sigma , b ]$ where $\sigma $ is an element of $\Sigma _n$ and 
$b$ is a list of $n$ binary digits representing the element of $(\mathbb{Z}_2)^n$. The multiplication then becomes 
$[\sigma , b ]\cdot [\tau  , c ]=[\sigma \cdot \tau , \tau ^{-1}(b)+c]$ where addition in the binary digits is a 
parity computation. We will use the convention of calling $b$ even or odd according to the number of 
ones appearing as binary digits of $b$. Notice that if $b$ and $c$ have the same parity, then $b+c$ is 
even, and if they differ in parity, then $b+c$ is odd.

Recall that the following sequence is a split exact sequence of groups:
$$(\mathbb{Z}_2)^n \stackrel{i}{\rightarrow} B_n \stackrel{\pi}{\rightarrow} \Sigma _n $$
where $i(b)=[1,b]$ and $\pi ([\sigma ,b])=\sigma $. 
We will also assume that $s: \Sigma _n \rightarrow B_n$ with $\pi \circ s= id_{\Sigma _n}$. 
Proper subgroups $G\subset B_n$ such that $\pi (G)=\Sigma _n$ are of particular interest here. Specifically we wish 
to examine maximal subgroups of $B_n$ with this property. Throughout this discussion, we will assume that 
$n\geq 5$. In order to examine the maximal subgroups of $B_n$ with $\pi (G)=\Sigma _n$, we need to introduce 
another maximal subgroup of $B_n$. Let $C_n$ be the subgroup of $B_n$ consisting of all elements $[\sigma ,b]$ 
such that $\sigma \in A_n$. It is immediately clear that $[B_n : C_n]=2$ and $C_n =\mathbb{Z}_2 \wr A_n$. So we get 
the following commutative diagram
$$\begin{CD} 
(\mathbb{Z}_2)^n @>{i }>> C_n @>{\pi }>> A_n \\
@|	                 @VVV		  @VVV\\
(\mathbb{Z}_2)^n @>{i }>> B_n @>{\pi }>> \Sigma _n \\
\end{CD}$$
where both horizontal sequences are split exact.

Let $G\subset B_n$ be a subgroup such that $\pi (G)=\Sigma _n$. Let $H=G \cap C_n$. We, thus, see that 
$\pi (H)=A_n$ and $[G:H]=2$. Letting $K=\ker (\pi |_H )$ gives the following commutative diagram:
$$\begin{CD} 
K @>>> H @>{\pi }>> A_n \\
@VVV	               @VVV              @| \\
(\mathbb{Z}_2)^n @>{i }>> C_n @>{\pi }>> A_n \\
\end{CD}$$
Notice that $C_n$ acts on $(\mathbb{Z}_2)^n$ by conjugation with kernel $i((\mathbb{Z}_2)^n)$. Since $K\lhd H$, 
we see that $K$ must be invariant under the action of $A_n$. Now we will look at the various possibilities for $K$ and 
the implications for $H$.

The first and most obvious case is that $K\cong (\mathbb{Z}_2)^n$ and $H=C_n$. In the other extreme 
lies the case where $K=1$, $H\cong A_n$. The third straightforward case is that $K=Z(C_n)$ 
and so $H\cong \mathbb{Z}_2\times A_n$. The fourth case is the most 
interesting one for this discussion and occurs when $K\cong (\mathbb{Z}_2)^{n-1}$ where each generator of $K$ is mapped 
to the product of two generators in $(\mathbb{Z}_2)^{n}\subset C_n$. Notice that $[C_n:H]=2$, so 
$H\lhd C_n$. Let $s$ be any splitting homomorphism of the split exact sequence 
defining $C_n$. Then $s(A_n)\cap H \lhd s(A_n)$ and $[s(A_n): H \cap s(A_n)]\leq 2$. Since $n\geq 5$, $A_n$ is simple, 
so $s(A_n)\subset H$. By realizing that these are the only possible cases for $K$, we have led ourselves to the 
following proposition:

\begin{proposition}\label{prop:subcn}
For $n\geq 5$, let $H$ be a subgroup of $C_n$ with $\pi (H)=A_n$; then $H$ is a split extension of $A_n$ by one of the 
following: $1$, $\mathbb{Z}_2$, $(\mathbb{Z}_2)^{n-1}$, or $(\mathbb{Z}_2)^n$. In the first two cases, $H$ is 
isomorphic to $A_n$ and $\mathbb{Z}_2 \times A_n$ respectively. In the last case, $H=C_n$, and in the third case 
$H=\{ [\sigma , b] \in C_n | b \textup{ is even} \}$.
\end{proposition}

We now use the following commutative diagram to relate Proposition \ref{prop:subcn} to the group $G$:
$$\begin{CD} 
K @>>> H @>{\pi }>> A_n \\
@|	               @VVV              @VVV \\
K @>>> G @>{\pi }>> \Sigma _n \\
@VVV	               @VVV              @| \\
(\mathbb{Z}_2)^n @>{i }>> B_n @>{\pi }>> \Sigma _n \\
\end{CD}.$$
Straightforward calculations using the fact that $[G:H]=2$ give us the following proposition:

\begin{proposition}\label{prop:subbn}
For $n\geq 5$, let $G$ be a subgroup of $B_n$ with $\pi (G)=\Sigma _n$; then $G$ is a split extension of $\Sigma _n$ by one of 
the following: $1$, $\mathbb{Z}_2$, $(\mathbb{Z}_2)^{n-1}$, or $(\mathbb{Z}_2)^n$. In the first two cases, $G$ is 
isomorphic to $\Sigma _n$ and $\mathbb{Z}_2 \times \Sigma _n$ respectively. In the last case $G=B_n$, and in the third 
case either $$G=\{ [\sigma , b] \in B_n | b \textup{ is even} \} \textup{  or  } G=\left\{ [\sigma , b] \in B_n \left| \begin{array}{c}
b \textup{ is even if }\sigma \in A_n\\
b \textup{ is odd if }\sigma \in \Sigma _n \setminus A_n
\end{array} \right. \right\}.$$
\end{proposition}

\section{Generators in $B_n$}

In order to find the strong symmetric genus of a group, it will be important to find a pair of generators. Of particular 
interest will be the orders of the generators in this pair as well as the order of their product.  In the following 
proposition, we will show that in the case of the hyperoctahedral groups, these orders will all be even.

\begin{proposition}
If $[\sigma , b]$ and $[\tau ,c]$ generate $B_n$, then the orders of $[\sigma , b]$, $[\tau ,c]$, and their product 
$[\sigma \cdot \tau, \tau ^{-1}(b)+c]$ all are even.
\end{proposition}

Proof: Suppose that one of the orders is odd. Without loss of generality, we may assume that $[\sigma , b]$ has 
odd order. Clearly $\sigma $ and $\tau $ must generate $\Sigma _n$. Since $[\sigma ,b]$ has odd order, 
$\sigma $ must also have odd order; therefore, $\sigma \in A_n$, $\tau \in \Sigma _n \setminus A_n$, and 
$\sigma \cdot \tau \in \Sigma _n \setminus A_n$. So we see that $\tau $ and $\sigma \cdot \tau$ have even orders.

Let $k$ be the order of $\sigma $. Since $[\sigma , b]$ has odd order, $[\sigma , b]^k=1$. Notice that $b$ must be 
even since $[\sigma , b]^k=[1,\sigma ^{-k+1} (b)+\sigma ^{-k+2}(b)+\cdots +\sigma ^{-1} (b)+b]$ and a sum of an odd 
number of odd elements is odd. If $c$ is also even, then $[\sigma , b]$ and $[\tau ,c]$ both lie in the proper 
subgroup $\{ [\sigma , b] \in B_n | b \textrm{ is even} \}\subset B_n$ and, therefore, cannot generate $B_n$. If on the 
other hand $c$ is odd, then $[\sigma , b]$ and $[\tau ,c]$ both lie in the proper subgroup 
$$\left\{ [\sigma , b] \in B_n \left| \begin{array}{c}
b \textrm{ is even if }\sigma \in A_n\\
b \textrm{ is odd if }\sigma \in \Sigma _n \setminus A_n
\end{array} \right. \right\}\subset B_n$$
and, therefore, cannot generate $B_n$. Thus the order of $[\sigma , b]$ must also be even. $\Box$

Notice that if $[\sigma , b]$ and $[\tau , c]$ generate $B_n$, then $\sigma $ and $\tau $ must generate $\Sigma _n$. 
Since we are looking to find generators of $B_n$, we need to first find generators of $\Sigma _n$. The 
second step in the process is to construct generators of $B_n$ from the generators of $\Sigma _n$. Also 
during this construction we would like to have control of the orders of the generators of $B_n$ as well 
as of their product. The following proposition gives the results of the process we will use to construct 
generators of $B_n$ from those of $\Sigma _n$.

\begin{proposition}\label{prop:genbn}
Suppose $\sigma $ and $\tau $ generate $\Sigma _n$ as $Symm(\Gamma )$ where $\Gamma =\{ 1,2,\dots ,n\}$ 
such that $\sigma \cdot \tau \in \Sigma _n \setminus A_n$ 
where $\sigma $, $\tau $, and $\sigma \cdot \tau $ all have even order. Futhermore, assume that $\sigma $ fixes an 
element $i\in \Gamma $ and $\tau $ fixes at least three elements of $\Gamma $, one of which we call $j$ such that $i$ and 
$j$ are in the same cycle of the element $\sigma \cdot \tau $. Let $b=(0,\dots , 0,1,0,\dots ,0)$ where the $1$ is in the 
$i^{th}$ position, and let 
$c=(0,\dots , 0,1,0,\dots ,0)$ where the $1$ is in the $j^{th}$ position. Under these conditions, $[\sigma , b]$ and 
$[\tau , c]$ generate $B_n$. In addition the elements $[\sigma , b]$, $[\tau , c]$ and $[\sigma \cdot \tau , \tau ^{-1} (b)+c]$
have the same orders as $\sigma $, $\tau $, and $\sigma \cdot \tau $ respectively.
\end{proposition}

Proof: The results about the orders of $[\sigma , b]$ and $[\tau , c]$ are obvious since $\sigma (b)=b$ and $\tau (c)=c$. 
On the other hand, $(\sigma \cdot \tau )^{-k}(\tau ^{-1}(b)+c)=(\sigma \cdot \tau )^{-k-1}(b)+(\sigma \cdot \tau )^{-k}(c)$. 
Notice that if $k$ is the length of the cycle in $\sigma \cdot \tau $ that contains $i$ and $j$, then 
$[\sigma \cdot \tau , \tau ^{-1}(b)+c ]^k = [(\sigma \cdot \tau)^k , (0,0, \dots ,0,0)]$ since the images of $b$ will put a 
$1$ in each position corresponding to the elements in the cycle and the images of $c$ will put a $1$ in each position 
corresponding to the elements in the same cycle. Thus it is clear that the order of $[\sigma \cdot \tau , \tau ^{-1}(b)+c ]$ 
is the same as the order of $\sigma \cdot \tau $.

Let $G= \langle [\sigma , b], [\tau , c] \rangle \subset B_n$. We need to show that $G=B_n$. Notice that $G$ is a subgroup 
of $B_n$ such that $\pi (G)=\Sigma _n$. So from Proposition \ref{prop:subbn} we know that $G$ is a split extension of 
$\Sigma _n$ by one of the following: $1$, $\mathbb{Z}_2$, $(\mathbb{Z}_2)^{n-1}$, or $(\mathbb{Z}_2)^n$.

Recall that any section $s:\Sigma _n \rightarrow B_n$ takes $\alpha \in A_n$ to $[1,d]\cdot 
[\alpha ,(0,\dots ,0)]\cdot [1,d] ^{-1}$ for some $[1,d] \in B_n$, and if $\alpha \in \Sigma _n \setminus A_n$, 
then either $s(\alpha )=[1,d]\cdot [\alpha ,(0,\dots ,0)]\cdot [1,d] ^{-1}$ or  $s(\alpha )=[1,d]\cdot [\alpha ,(1,\dots ,1)]
\cdot [1,d] ^{-1}$. Notice if $\alpha \in A_n$, $s(\alpha )=[\alpha ,a]$ where $a$ is even.

Suppose first that $G$ is a split extension of $\Sigma _n$ by $1$. Then $[\sigma , b]$ and $[\tau , c]$ must be 
in the image of some section homorphism $s:\Sigma _n \rightarrow B_n$. Either $\sigma \in A_n$ 
or $\tau \in A_n$, however, and both $b$ and $c$ are odd. So $G$ is not a split extension of $\Sigma _n$ by $1$.

Suppose that $G$ is a split extension of $\Sigma _n$ by $\mathbb{Z}_2$. This implies that for some 
section homorphism, $s:\Sigma _n \rightarrow B_n$, $G=Z(B_n) \times s(\Sigma _n)$. If $n$ is even, then 
the only elements in $Z(B_n) \times s(\Sigma _n)$ of the form $[\alpha ,a]$ with $\alpha \in A_n$ have an even 
$a$. Clearly either $[\sigma , b]\not \in Z(B_n) \times s(\Sigma _n)$ or $[\tau , c]\not \in Z(B_n) \times 
s(\Sigma _n)$; therefore, we may assume that $n$ is odd.  Suppose that $\tau \in A_n$. So $[\tau , c]=
[1,(1,\dots ,1)]\cdot [1,d]\cdot [\tau ,(0,\dots ,0)] \cdot [1,d] ^{-1}$ for some $[1,d] \in B_n$. It follows that 
$c=(1,\dots ,1) + \tau ^{-1} (d)+d $, which cannot be the case since $\tau $ fixes an element besides $j$. 
So we may also assume that $\tau \in \Sigma _n \setminus A_n $.

Since $\tau \in \Sigma _n \setminus A_n $, either $[\tau , c]=[1,(1,\dots ,1)]\cdot [1,d]\cdot [\tau ,(0,\dots ,0)] 
\cdot [1,d] ^{-1}$ or $[\tau , c]= [1,d]\cdot [\tau ,(1,\dots ,1)] \cdot [1,d] ^{-1}$ for some $[1,d] \in B_n$. 
In either case, $c=(1,\dots ,1) + \tau ^{-1} (d)+d $, which as before cannot be the case since $\tau $ fixes 
an element besides $j$. So $G$ is not a split extension of $\Sigma _n$ by $\mathbb{Z}_2$.

Suppose that $G$ is a split extension of $\Sigma _n$ by $(\mathbb{Z}_2)^{n-1}$. By Proposition 
\ref{prop:subbn}, either 
$$G=\{ [\sigma , b] \in B_n | b \textrm{ is even} \} \textrm{  or  } G=\left\{ [\sigma , b] \in B_n \left| \begin{array}{c}
b \textrm{ is even if }\sigma \in A_n\\
b \textrm{ is odd if }\sigma \in \Sigma _n \setminus A_n
\end{array} \right. \right\}.$$
Recall that either $\sigma \in A_n$ or $\tau \in A_n$. If $\sigma \in A_n$, then $[\sigma ,b] \not \in G$. If $\tau \in A_n$, 
then $[\tau ,c] \not \in G$. So $G$ is not a split extension of $\Sigma _n$ by $(\mathbb{Z}_2)^{n-1}$.

The only possibility left is that $G$ is a split extension of $\Sigma _n$ by $(\mathbb{Z}_2)^{n}$. Thus we see that 
$G=B_n$ and so the proposition is proven. $\Box $

\section{generating pairs}
If a finite group $G$ has generators $x$ and $y$ of orders $p$ and $q$ respectively with $xy$ having the 
order $r$, then we say that $(x,y)$ is a $(p,q,r)$ generating pair of $G$. By the obvious symmetries concerning 
generators, we will use the convention that $p\leq q\leq r$. Following the convention of Marston Conder in 
\cite{c:srqtg}, we say that a $(p,q,r)$ generating pair of $G$ is a minimal generating pair if there does not 
exist a $(k,l,m)$ generating pair for $G$ with $\frac{1}{k}+\frac{1}{l}+\frac{1}{m} >\frac{1}{p} +\frac{1}{q} +
\frac{1}{r}$. 

Recall that the groups of small strong symmetric genus are well known 
(see \cite{b:cfga,mz:gsssg}). So for now we will assume that $\sigma ^0 (G) >1$. The hyperoctahedral groups 
have strong symmetric genus at least 2. It is also known that for groups with $\sigma ^0 (G) >1$, any 
generating pair will be a $(p,q,r)$ generating pair with $\frac{1}{p} +\frac{1}{q} + \frac{1}{r}<1$. 
Using the Riemann-Hurwitz equation, we see that given any generating pair of $G$, we get an upper bound 
on the strong symmetric genus of $G$ \cite{s:sfg}. If $G$ has a $(p,q,r)$ generating pair, then 
$$\sigma ^0 (G)\leq 1+\frac{1}{2} |G|\cdot (1-\frac{1}{p} -\frac{1}{q} -\frac{1}{r}).$$ The following lemma, which is a result of 
Singerman \cite{s:srs} (see also \cite{mz:gsssg,t:fgas}), shows that the strong symmetric genus for most groups 
is computed directly from a minimal generating pair.

\begin{lemma}[Singerman \cite{s:srs}]\label{lemma:sing}
Let $G$ be a finite group such that $\sigma ^0 (G)\geq 2$. If $|G|>12(\sigma ^0 (G) -1)$, then $G$ has a 
$(p,q,r)$ generating pair with $$\sigma ^0 (G)=1+\frac{1}{2} |G|\cdot (1-\frac{1}{p} -\frac{1}{q} -\frac{1}{r}).$$ 
In addition, we may assume that either $p=2$, or $p=q=3$ and $r$ is $4$ or $5$.
\end{lemma}

Applying Lemma \ref{lemma:sing} to our situation with the hyperoctahedral groups results in the following corollary.

\begin{corollary}\label{cor:sing}
Let $n\geq 3$ and $G=B_n$. If $|G|>12(\sigma ^0 (G) -1)$, then $G$ has a  
$(p,q,r)$ generating pair with $$\sigma ^0 (G)=1+\frac{1}{2} |G|\cdot (1-\frac{1}{p} -\frac{1}{q} -\frac{1}{r}).$$ 
In addition, we may assume that $p=2$, $q=4$, and $r$ is either 6, 8 or 10.
\end{corollary}

Proof: Since all groups $H$ with $\sigma ^0 (H) < 2$ have been classified, we know that $\sigma ^0 (G)\geq 2$. 
We also know that if $G$ has a $(p,q,r)$ generating pair, then $p$, $q$, and $r$ must all be even. 
The rest of the corollary follows from straightforward calculations and by recalling that $|G|=n! 2^n$. $\Box$

In light of Corollary \ref{cor:sing}, Corollary \ref{cor:main} is a direct result of Theorem \ref{thm:main}. 
We will see that the exceptional cases of Theorem \ref{thm:main} will be shown by using GAP \cite{g:gap}. 
For the remaining cases, we will show the existence of 
a $(2,4,6)$ generating pair for $B_n$ by constructing it from a $(2,4,6)$ generating pair for $\Sigma _n$ using 
Proposition \ref{prop:genbn}.

\section{special generators of $\Sigma _n$}
First we will need to recall a theorem of Jordan (see \cite{w:fpg}).

\begin{theorem}[Jordan \cite{w:fpg}]\label{thm:jordan}
Let $G$ be a transitive permutation group on a set of size $n$ such that $G$ is primitive and contains a 
$p$-cycle for some prime $p<n-2$. Then $G$ is either $A_n$ or $\Sigma _n$.
\end{theorem}

We will use this theorem by applying the following corollary. The method that we will use follows that of 
Marston Conder (see \cite{c:gasg,c:mgasg}).

\begin{corollary}\label{cor:gensig}
Let $G=\langle \sigma , \tau  \rangle $ be a transitive permutation group on a set $\Gamma $ of size $n$. Assume 
that for some prime $p<n-2$, $G$ contains a $p$-cycle $z$ such that the points permuted by $z$ include an 
element $a\in \Gamma $ as well as 
its image under $\sigma $ and an element $b\in \Gamma $ as well as its image under $\tau $. Then $G$ is either 
$A_n$ or $\Sigma _n$.
\end{corollary}

Proof: Using Theorem \ref{thm:jordan}, the only part we need to prove is that $G$ is a primitive permutation 
group on $\Gamma $. Suppose that $G$ acts imprimitively on $\Gamma $. Let $z$ be the $p$-cycle 
in the hypothesis. Since $z$ fixes all the other points of $\Gamma $, the $p$ points permuted by $z$ 
must lie in the same block of imprimitivity $B$. By hypothesis, however, $a$ and $\sigma (a)$ are both in 
$B$; therefore, $\sigma $ fixes $B$. Also by hypothesis, $b$ and $\tau (b)$ are both in $B$; therefore, 
$\tau $ fixes $B$. This implies that the block $B$ is fixed by all of $G$, contradicting the transitivity of $G$. 
$\Box $

At this point we have done all of the preliminary work; it is now time to prove the main theorem. To begin we 
point out that the result for the groups $B_n$ with $3\leq n \leq 40$ were shown using GAP \cite{g:gap}. Some of 
these results were shown in \cite{vwg:ggt}.
For the groups $B_n$ with $n \geq 41$, we show that each group $\Sigma _n$ with $n \geq 41$ can be 
generated by two elements $\sigma $ and $\tau $. Using Corollary \ref{cor:gensig}, this process will involve 
finding two elements $\sigma $ and $\tau $ such that these two elements generate a transitive permutation group on $n$ 
letters such that some element of $\langle \sigma ,\tau \rangle$ is a $p$-cycle, which contains a letter 
and its image under $\sigma $ and a possibly different letter and its image under $\tau $. In order to 
be able to construct a (2,4,6) generating pair, $([\sigma ,b],[\tau ,c])$, of $B_n$ from 
$\sigma $ and $\tau $ as in Proposition \ref{prop:genbn}, we will require that the following hold:
\begin{enumerate}
\item $\sigma $ has order 4 and $\tau $ has order 2,
\item $\sigma $ fixes a letter $i$,
\item $\tau $ fixes at least three letters one of which we call $j$, 
\item $\sigma \cdot \tau \in \Sigma _n \setminus A_n $, 
\item $i$ and $j$ are in the same cycle of the element $\sigma \cdot \tau $,
\item $\langle \sigma ,\tau \rangle$ is a transitive permutation group on $n$ letters, and 
\item  some element of $\langle \sigma ,\tau \rangle$ is a $p$-cycle, which contains a letter 
and its image under $\sigma $ and a possibly different letter and its image under $\tau $.
\end{enumerate}

\section{Coset Diagrams and Results}

\begin{figure}
\begin{center}
\includegraphics{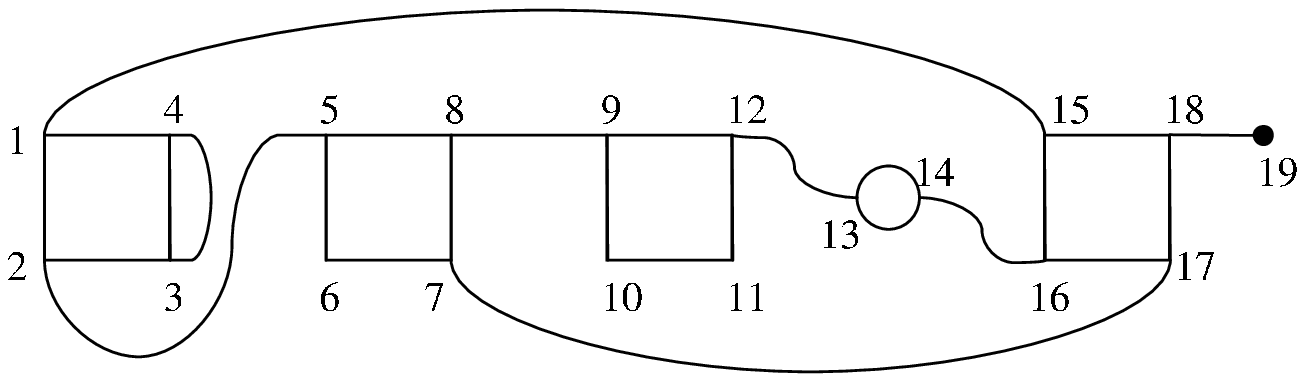}
\caption{Example $\Sigma _{19}$}
\end{center}
\end{figure}

We will demonstrate the existence of the desired elements $\sigma $ and $\tau $ in each $\Sigma _n$ 
with $n\geq 41$ by using diagrams which are called coset diagrams. First let us give a brief introduction to these 
diagrams by looking at an example in Figure 1. The diagram in Figure 1 is a representation of the group 
$\Sigma _{19}$. The diagram gives two elements of a transitive permutation group on 19 letters. 
The first element $\sigma $ is given by the squares and circles in the diagram and will have order 4. A cycle 
of $\sigma $ is given by reading the labels in a counter-clockwise direction. So in this case, 
$$\sigma =(1\;2\;3\;4)(5\;6\;7\;8)(9\;10\;11\;12)(13\;14)(15\;16\;17\;18).$$
The second element $\tau $ is given by the other lines and curves in the diagram. $\tau $ 
will be an involution with a 2-cycle, for each line or curve, that is the cycle of its two labels. In this case we see
$$\tau =(1\;15)(2\;5)(3\;4)(7\;17)(8\;9)(12\;13)(14\;16)(18\;19).$$
To find to where an integer is mapped under the element $\sigma \cdot \tau $, we start at that label, trace along 
a square or circle in the counter-clockwise direction if possible, and then trace a curve not in one of these shapes 
if possible. The permutation given by a general word in $\sigma $ and $\tau $ can be found in a similar way.

We see that the element $\sigma ^2 \cdot \tau $ has four 4-cycles and one 3-cycle, which is $(10\;13\;12)$. By Corollary 
\ref{cor:gensig}, we see that $\sigma $ and $\tau $ generate the symmetric group on 19 elements. 
Using the elements $\sigma $ and $\tau $ and creating elements $[\sigma ,b]$ and $[\tau ,c]$ of $B_{19}$ as in 
Proposition \ref{prop:genbn}, we also see that $[\sigma ,b]$ and $[\tau ,c]$ generate $B_{19}$.

\begin{figure}
\begin{center}
\includegraphics{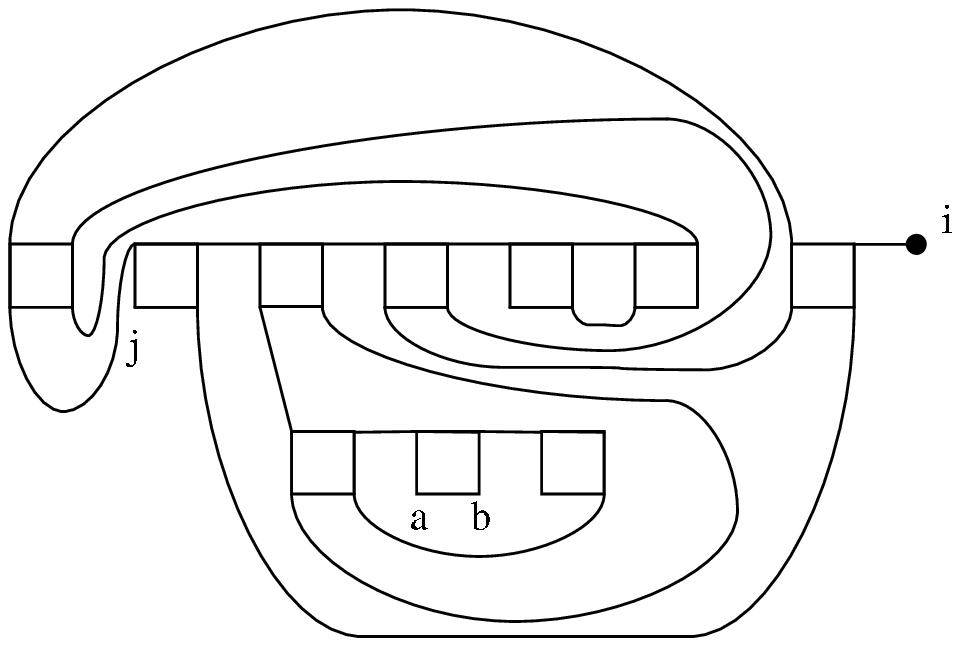}
\caption{$\Sigma _{n}$ for $n\equiv 41 (12)$ and $n\geq 41$}
\end{center}
\end{figure}

\begin{figure}
\begin{center}
\includegraphics{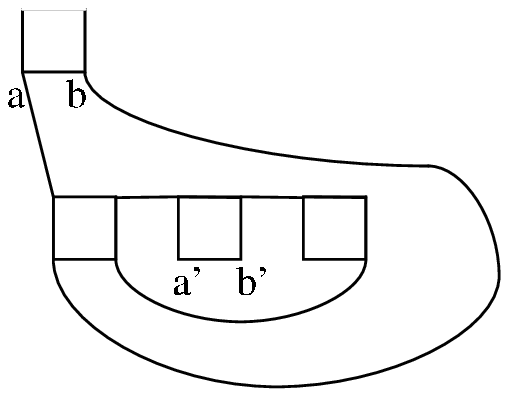}
\caption{additional piece}
\end{center}
\end{figure}

In a similar way as the example above, we will now look at Figure 2, which gives a coset digram for the group 
$\Sigma _{41}$. As with the example in Figure 1, these figures show that the group 
$B_{41}$ is a $(2,4,6)$ group. The necessary information about $\langle \sigma ,\tau \rangle$ 
containing a $p$-cycle, which contains a letter and its image under $\sigma $ and a possibly different letter and its 
image under $\tau $, can be found in Table I. Similar coset diagrams exist for the groups 
$\Sigma _n$ with $42\leq n \leq 52$ and the corresponding information is also contained in Table I. These diagrams 
show that each group $B_n$ for $42\leq n \leq 52$ is a $(2,4,6)$ groups.

\begin{table}
\begin{tabular}{|c|c|c|c|c|}\hline \label{tab:one}
$n$ value & \multicolumn{2}{c}{conditions on $m$} \hfill \vline & word in $\sigma $ and $\tau $ & prime \\
\hline
\hline
$n=41+12m$ & $m\not \equiv 0 (13)$ & & $\sigma \tau \sigma \tau \sigma $ & 13 \\
\hline
$n=41+12m$ & $m\not \equiv 8 (13)$ & $m\not \equiv 9 (13)$ & $\sigma ^3  \tau  \sigma  \tau  \sigma 
 \tau  \sigma ^2 \tau $ & 13 \\
\hline
$n=42+12m$ & $m\not \equiv 1 (7)$ & & $\sigma  \tau  \sigma ^3  \tau $ & 7 \\
\hline
$n=42+12m$ & $m\not \equiv 2 (7)$ & $m\not \equiv 4 (7)$ & $\sigma ^3  \tau  \sigma  \tau  \sigma ^2
 \tau  \sigma  \tau  \sigma ^2 \tau \sigma ^3 \tau  \sigma ^2$ & 7 \\
\hline
$n=43+12m$ & $m\not \equiv 4 (7)$ & & $\tau  \sigma  \tau  \sigma ^3$ & 7 \\
\hline
$n=43+12m$ & $m\not \equiv 2 (7)$ & $m\not \equiv 3 (7)$ & $\tau  \sigma ^2  \tau  \sigma  \tau  
\sigma  \tau  \sigma ^2 \tau  \sigma  \tau  \sigma ^2 \tau  \sigma  \tau $ & 7 \\
\hline
$n=44+12m$ & $m\not \equiv 1 (13)$ & & $\tau  \sigma  \tau  \sigma ^3 \tau  \sigma  \tau  
\sigma ^3 \tau  \sigma $ & 13 \\
\hline
$n=44+12m$ & $m\not \equiv 9 (13)$ & $m\geq 1$ & $\sigma ^2  \tau  \sigma ^2 \tau  \sigma  \tau  
\sigma ^3 \tau  \sigma  \tau  \sigma  \tau  \sigma ^3$ & 13 \\
\hline
$n=45+12m$ & $m\not \equiv 5 (7)$ & $m\not \equiv 6 (7)$ & $\sigma ^2 \tau \sigma ^2\tau \sigma $ & 7 \\
\hline
$n=45+12m$ & $m\not \equiv 1 (7)$ & $m\not \equiv 3 (7)$ & $\tau \sigma ^3  \tau  \sigma  \tau  \sigma 
 \tau  \sigma ^3 \tau \sigma ^2$ & 7 \\
\hline
$n=46+12m$ & $m\not \equiv 0 (7)$ & & $\tau \sigma ^3\tau \sigma ^2\tau \sigma \tau \sigma ^2\tau \sigma ^3
\tau \sigma ^2\tau \sigma \tau \sigma ^3\tau \sigma \tau $ & 7 \\
\hline
$n=46+12m$ & $m\not \equiv 1 (7)$ & $m\not \equiv 5 (7)$ & $\tau \sigma ^2  \tau  \sigma ^3\tau  \sigma 
 \tau  \sigma \tau \sigma ^3\tau \sigma $ & 7 \\
\hline
$n=47+12m$ & $m\not \equiv 5 (7)$ & $m\not \equiv 6 (7)$ & $\sigma ^2 \tau \sigma \tau \sigma \tau \sigma ^3
\tau $ & 7 \\
\hline
$n=47+12m$ & $m\not \equiv 1 (7)$ & $m\not \equiv 3 (7)$ & $\sigma \tau \sigma ^2\tau  \sigma  \tau  \sigma 
 \tau  \sigma ^2 \tau \sigma ^3\tau $ & 7 \\
\hline
$n=48+12m$ & $m\not \equiv 1 (13)$ & & $\sigma \tau \sigma \tau \sigma $ & 13 \\
\hline
$n=48+12m$ & $m\not \equiv 11 (13)$ & $m\geq 1$ & $\sigma \tau \sigma \tau \sigma ^2 
\tau  \sigma ^3 \tau $ & 13 \\
\hline
$n=49+12m$ & $m\not \equiv 5 (7)$ &$m\not \equiv 6 (7)$ & $\sigma ^2\tau \sigma \tau \sigma \tau \sigma ^3
\tau $ & 7 \\
\hline
$n=49+12m$ & $m\not \equiv 1 (7)$ & $m\not \equiv 3 (7)$ & $\sigma ^2  \tau  \sigma  \tau  \sigma ^3
\tau  \sigma \tau \sigma ^3 \tau \sigma \tau \sigma ^2\tau \sigma $ & 7 \\
\hline
$n=50+12m$ & $m\not \equiv 4 (7)$ & & $\tau \sigma \tau \sigma \tau \sigma ^2\tau \sigma ^2\tau \sigma 
\tau \sigma \tau $ & 7 \\
\hline
$n=50+12m$ & $m\not \equiv 6 (7)$ & $m\geq 2$ & $\sigma ^2\tau \sigma \tau \sigma ^3 \tau \sigma \tau 
\sigma ^2\tau  \sigma ^3 \tau \sigma ^3$ & 7 \\
\hline
$n=51+12m$ & $m\not \equiv 6 (7)$ & & $\sigma ^3\tau \sigma ^2\tau \sigma \tau \sigma \tau \sigma ^2
\tau \sigma ^2\tau $ & 7 \\
\hline
$n=51+12m$ & $m\not \equiv 3 (7)$ & & $\sigma \tau \sigma ^2\tau  \sigma \tau \sigma ^2 \tau \sigma \tau 
\sigma ^3\tau  \sigma \tau \sigma ^3 $ & 7 \\
\hline
$n=52+12m$ & $m\not \equiv 6 (11)$ & & $\tau \sigma ^2\tau \sigma \tau \sigma ^2\tau \sigma ^3$ & 11 \\
\hline
$n=52+12m$ & $m\not \equiv 10 (11)$ & $m\geq 1$ & $\sigma ^2\tau \sigma ^2\tau  \sigma \tau \sigma \tau 
\sigma ^2$ & 11 \\
\hline
\end{tabular}
\caption{Demonstration of $p$-cycles}
\end{table}

To show the result that the groups $B_n$ with $n > 52$ are all $(2,4,6)$ groups, we will recursively add the piece shown 
in Figure 3. We attach the piece by matching up the points $a$ and $b$ and then making the $a'$ and $b'$ into a 
new $a$ and $b$ respectively. Attaching this piece recursively to the diagram in Figure 2 gives a diagram 
for each $n \geq 41$ with $n \equiv 41 (\textrm{mod } 12)$. Again the 
necessary information about $\langle \sigma ,\tau \rangle$ containing a $p$-cycle, which contains a letter and its 
image under $\sigma $ and a possibly different letter and its image under $\tau $, can be found in Table I. 
Each $n\geq 41$ with $n \equiv 41  (\textrm{mod } 12)$ satisfies the conditions of at least one of the rows of Table I. 
The word in $\sigma $ and $\tau $ in that row then has some power which is a $p$-cycle for the prime $p\;$ in 
the last column. A similar recursive procedure on the diagrams for the groups $\Sigma _n$ with 
$42\leq n \leq 52$ together with the information in Table I shows that the remaining groups $B_n$ with 
$n\geq 54$ are also $(2,4,6)$ groups. This completes the demonstration that 
every $B_n$ with $n\geq 41$ is a $(2,4,6)$ group and finishes the proof of the main theorem. $\Box$

\end{document}